\documentclass[twoside]{article}
\usepackage{graphics}

\author{Oleg Pikhurko\thanks{Supported by a Senior Rouse Ball
 Studentship, Trinity College, Cambridge, UK.}\\
 DPMMS, Centre for Mathematical Sciences\\
 Cambridge University, Cambridge~CB3~0WB, England\\
 E-mail: {\tt O.Pikhurko@dpmms.cam.ac.uk}}

\usepackage{amssymb}
\newcommand{\eqref}[1]{\mbox{\rm(\ref{#1})}}
\newcommand{\C}[1]{{\protect\cal #1}}

\newcommand{\I}[1]{{\mathbb #1}}

\newcommand{\binom}[2]{{#1\choose #2}}
\newcommand{\e}{\epsilon}
\newcommand{\me}{\mathrm{e}}
\newcommand{\md}{\mathrm{d}}

\newcommand{\qed}{\nolinebreak\mbox{\hspace{5 true pt}%
\rule[-0.85 true pt]{3.9 true pt}{8.1 true pt}}}

\newtheorem{theorem}{Theorem}
\newtheorem{lemma}[theorem]{Lemma}

\newtheorem{problem}[theorem]{Problem}

\newcommand{\ceil}[1]{\lceil #1\rceil}
\newcommand{\floor}[1]{\lfloor #1\rfloor}
\begin{document}

\newcommand{\magic}{\C M}
\newcommand{\imagic}{\C I}
\newcommand{\amagic}{\C A}
\newcommand{\MagicUpper}{0.489...}
\newcommand{\BSupUpper}{0.489...}
\newcommand{\Q}[1]{{\tt IF[\mbox{#1}]}}

\renewcommand{\arraystretch}{1.3}
\renewcommand{\baselinestretch}{1.25}\rm

\title{Dense Edge-Magic Graphs and\\ Thin Additive Bases} 
\author{Oleg Pikhurko%
 \footnote{A part of this research was carried out when the author was
   supported by a Research Fellowship of St.~John's College,
   Cambridge, UK.}\\
 Department of Mathematical Sciences\\
 Carnegie Mellon University\\
 Pittsburgh, PA 15213-3890}
 \maketitle

\pagestyle{myheadings}
\markboth{\sc O.~Pikhurko}{\sc Edge-Magic Graphs and Additive Bases}

\begin{abstract}
 A graph $G$ of order $n$ and size $m$ is {\em edge-magic} if there is
a bijection $l:V(G)\cup E(G)\to [n+m]$ such that all sums
$l(a)+l(b)+l(ab)$, $ab\in E(G)$, are the same. We present new lower
and upper bounds on $\magic(n)$, the maximum size of an edge-magic
graph of order $n$, being the first to show an upper bound of the form
$\magic(n)\le (1-\e)\binom n2$. Concrete estimates for $\e$ can be
obtained by knowing $s(k,n)$, the maximum number of distinct pairwise
sums that a $k$-subset of $[n]$ can have.

So, we also study $s(k,n)$, motivated by the above connections to
edge-magic graphs and by the fact that a few known functions from
additive number theory can be expressed via $s(k,n)$. For example, our
estimate
 $$
 s(k,n)\le n+k^2\left(\frac14-\frac1{(\pi+2)^2}+o(1)\right)
 $$
 implies new bounds on the maximum size of quasi-Sidon
sets, a problem posed by Erd\H os and Freud [J.\ Number
Th. \textbf{38} (1991) 196--205]. The related problem for differences
is considered as well.

\vspace{10pt}\noindent

\vspace{10pt}\noindent
 \end{abstract}

\section{Introduction}

Let $[k]$ stand for $\{1,\dots,k\}$.  Let $G$ be a graph with $n$
vertices and $m$ edges.  An {\em edge-magic labelling} with the
{\em magic sum $s$} is a bijection $l:V(G)\cup E(G)\to [m+n]$ such that
$l(a)+l(b)+l(ab)=s$ for any edge $ab$ of $G$. (We always assume that
$V(G)\cap E(G)=\emptyset$.) This definition appeared first in Kotzig
and Rosa~\cite{kotzig+rosa:70} (but under the name {\em magic
valuation}). The graph $G$ is {\em edge-magic} if it admits an
edge-magic labelling (for some $s$). We refer the reader to
Gallian~\cite{gallian:DS} and Wood~\cite{wood:02} for plentiful
references on edge-magic graphs.

Not all graphs are edge-magic nor is this property in any way monotone
with respect to the subgraph relation. In 1996 Erd\H os
 asked (see~\cite{craft+tesar:99}) for $\magic(n)$, the maximum
number of edges that an edge-magic graph of order $n$ can have.

This function has been computed exactly for $n\le 6$ but for large $n$
the best known bounds were $\floor{n^2/4}\le \magic(n)\le \binom n2
-1$, see Craft and Tesar~\cite{craft+tesar:99}.

Here we improve both these bounds if $n$ is large.

\begin{theorem}\label{th:magic}
 \begin{equation}\label{eq:magic}
 \frac27\, n^2+O(n)\le \magic(n)\le
\left(\MagicUpper+o(1)\right)\, n^2.
 \end{equation}
 \end{theorem}

It turns out that edge-magic labellings have strong relations to some
problems from additive number theory, especially to additive bases.

Section~\ref{Wood} can serve as a warm-up where we improve the bounds
of Wood~\cite{wood:02} on so-called {\em edge-magic injections}. Our
proof uses some classical results about {\em Sidon sets}, that is, sets
$A\subset \I Z$ such that all sums $a+b$, with $a,b\in A$ and $a\le
b$, are distinct.

For a set $A$ of integers define its {\em sum-set} by $A+A:=\{a+b\mid
a,b\in A\}$; $A$ is called an {\em additive basis} for $X$ if
$A+A\supset X$.  In Section~\ref{lower}, we prove the lower bound
in~\eqref{eq:magic} by using known (explicit) constructions of a thin
additive basis for some suitable interval of integers.

But the most interesting connections were found during our quest for
an upper bound on $\magic(n)$. This research led us to the following
problem. What is
 $$
 s(k,n):=\max\left\{|A+A|\mid A\in \binom{[n]}k\right\},
 $$
 that is, the maximum size of the sum-set of a $k$-subset of
$\{1,\dots,n\}$?

The trivial upper bound is 
 \begin{equation}\label{eq:trivial}
 s(k,n)\le \min\left\{\binom k2+k,2n-1\right\}.
 \end{equation}

We have $s(k,n)=\binom k2+k$ if and only if there exists a Sidon
$k$-set $A\subset [n]$; the classical results of
Singer~\cite{singer:38} and Erd\H os and Tur\'an~\cite{erdos+turan:41}
(see~\cite[Chapter~II]{halberstam+roth:s}) state that for a given $n$
the largest such $k$ is $(1+o(1))\, n^{1/2}$. The open question
whether the maximum size of a Sidon subset of $[n]$ is $n^{1/2}+O(1)$
has the \$500-dollar reward of Erd\H os~\cite{erdos:80} attached.

We have $s(k,n)=2n-1$ if and only if there is an additive $k$-basis
$A\subset[n]$ for $[2,2n]$. How small can $k$ be then? A simple
construction of Rohrbach~\cite[Satz~2]{rohrbach:37} gives $(2\sqrt
2+o(1))\, n^{1/2}$ for $k$ (see Section~\ref{SLower}). The trivial
lower bound is $k\ge (2+o(1))\, n^{1/2}$; the current best known bound
$k\ge (2.17...+o(1))\, n^{1/2}$ of Moser, Pounder and
Riddell~\cite{moser+pounder+riddell:69} is only slightly bigger.

As we see, already the question when we have equality in~\eqref{eq:trivial}
leads to very difficult open problems. The computation of $s(k,n)$ for
other values is likely to be even harder. We present the following
upper bound which improves on~\eqref{eq:trivial} for a range of $k$ around
$2\, n^{1/2}$. 

\begin{theorem}\label{th:SUpper} 
 \begin{equation}\label{eq:SUpper}
 s(k,n)\le n+k^2\left(\frac14-\frac1{(\pi+2)^2}+o(1)\right).
 \end{equation}\end{theorem}

Here is an application of Theorem~\ref{th:SUpper}. Erd\H os and
Freud~\cite{erdos+freud:91} call a set $A\in\binom{[n]}k$ with
$|A+A|=(1+o(1))\binom k2$ {\em quasi-Sidon} and ask how large $k$ can
be. (It is obvious what is meant here so we do not bother writing out
any formal definitions.)  They constructed quasi-Sidon subsets of
$[n]$ with
 \begin{equation}\label{eq:QuasiSidonConstruction}
 k=(2/\sqrt3+o(1))\,n^{1/2}=(1.154...+o(1))\,n^{1/2}.
 \end{equation}
 As $A+A\subset[2n]$, a trivial upper bound is $\binom k2\le
(2+o(1))\,n$, that is, $k\le (2+o(1))\,n^{1/2}$. Erd\H os and
Freud~\cite[p.~204]{erdos+freud:91} promised to publish the proof of
$k\le (1.98+o(1))\,n^{1/2}$ in a follow-up paper. Unfortunately, it
has not been published, but their bound is superseded by the following
easy corollary of Theorem~\ref{th:SUpper} anyway.

\begin{theorem}\label{th:QuasiSidon} Let $A\subset [n]$ be quasi-Sidon. Then
 $$
 |A|\le 
\left(\left(\frac14+\frac{1}{(\pi+2)^2}\right)^{-1/2}+o(1)\right)\, n^{1/2}=
(1.863...+o(1))\, n^{1/2}.\qed
 $$
 \end{theorem}

As another application of Theorem~\ref{th:SUpper} let us show that
$\magic(n)\le (1-\e) \binom n2$. Indeed, if $G$ is an edge-magic graph
of order $n$ and size $(\frac12+o(1))\, n^2$, then its vertex labels
form a quasi-Sidon set, which contradicts
Theorem~\ref{th:QuasiSidon}. This way we do not obtain any explicit value
for $\e$ but one can get one by using Theorem~\ref{th:SUpper} with a
little bit of work. A slightly better bound, the one in~\eqref{eq:magic},
is deduced in Section~\ref{upper} from a generalisation of
Theorem~\ref{th:SUpper}.

\begin{figure}[h]
\begin{minipage}[t]{6.0cm}
\begin{center}
\scalebox{0.59}[0.8]{\includegraphics{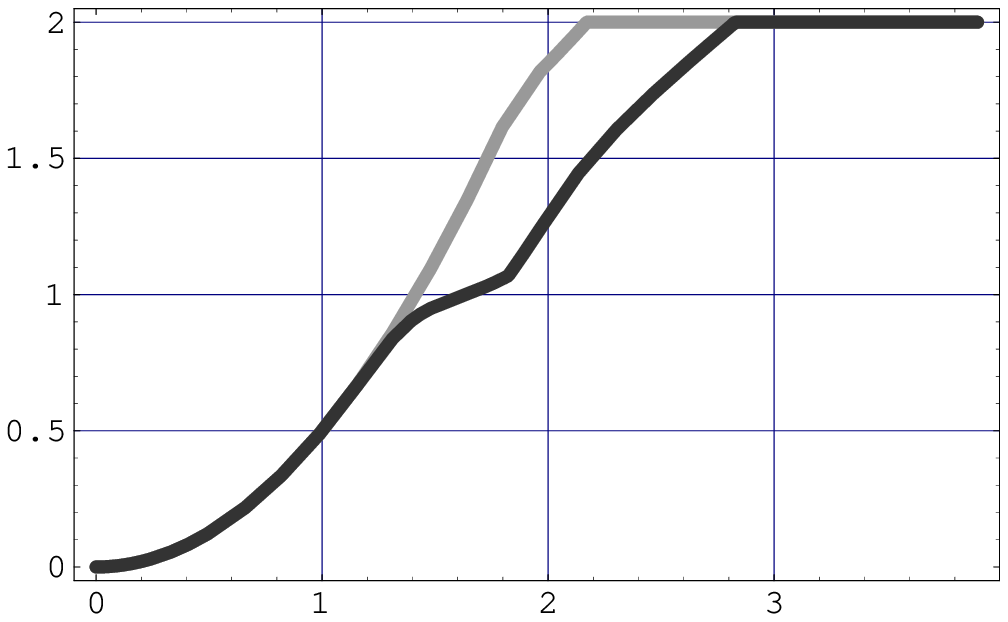}}
\end{center}
\end{minipage}
\hfill
\begin{minipage}[t]{6.0cm}
\begin{center}
\scalebox{0.59}[0.8]{\includegraphics{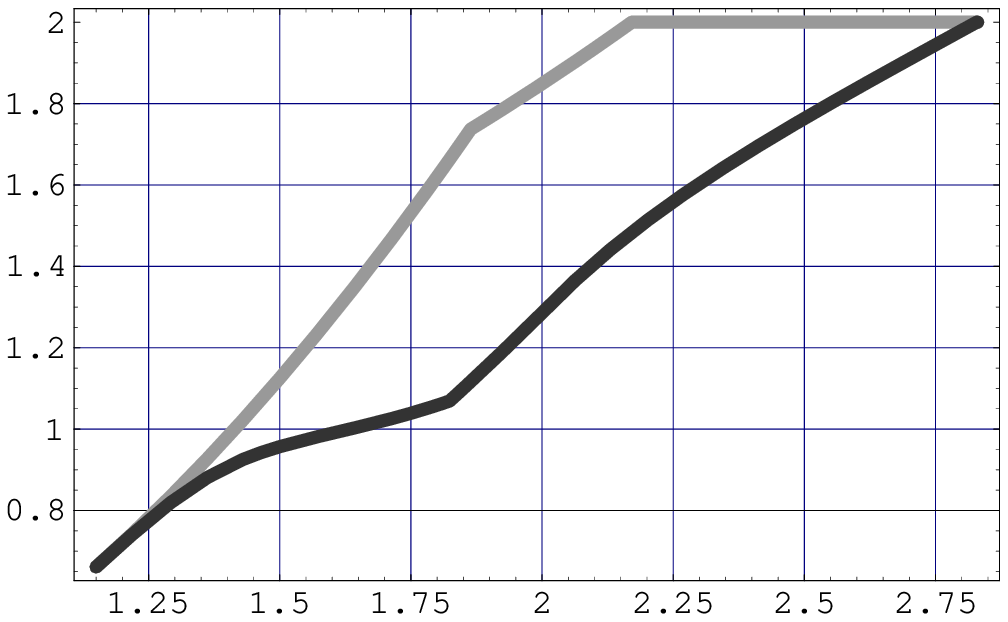}}
\end{center}
\end{minipage}
\caption{\small Our bounds on $s(k,n)$: $x=kn^{-1/2}$, 
 $y=s(k,n)/n$.}\label{fg:s}
\end{figure}

Given these applications of $s(k,n)$, we present some lower bounds on
$s(k,n)$ in Section~\ref{SLower}. It is interesting to compare them
with the upper bounds, see Figure~\ref{fg:s}.

Our auxiliary Lemma~\ref{lm:EF+L} states that any asymptotically maximum
Sidon subset of $[n]$ is uniformly distributed in subintervals and in
residue classes simultaneously. This places the corresponding results
of Erd\H os and Freud~\cite{erdos+freud:91} and
Lindstr\"om~\cite{lindstrom:98} under a common roof.

Besides being a natural and interesting question on its own, the
$s(k,n)$-problem demonstrates new connections between Sidon sets and
additive bases. This helped the author to realise that the technique
of Moser~\cite{moser:60} which was used in the context of additive
bases can be applied to $s(k,n)$ (and to quasi-Sidon sets). In fact,
our proof of Theorem~\ref{th:SUpper} goes by modifying
Moser's~\cite{moser:60} method. Although the determination of $s(k,n)$
is apparently very hard, it seems a promising direction of research.

In Section~\ref{differences} we study the analogous problem for
differences.

\section{Edge-Magic Injections}\label{Wood}

Wood~\cite{wood:02} defines an {\em edge-magic injection} of a graph
$G$ as an injection $l:V(G)\cup E(G)\to \I Z_{>0}$ (into positive
integers) such that for any edge $ab\in E(G)$ the sum
$l(a)+l(b)+l(ab)=s$ is constant. Note that the labels need not sweep a
contiguous interval of integers (but must be pairwise distinct). It is
easy to show that any graph $G$ admits an edge-magic injection.

The general question is how economical such a labelling can be.  One
possible way to state it formally is to ask about $\imagic(G)$, the
smallest value of the magic sum $s$ over all edge-magic injections of
$G$. If $v(G)=n$, then clearly $\imagic(G)\le \imagic(K_n)$, so here
we investigate $\imagic(K_n)$. Wood~\cite[Theorem~1]{wood:02} showed
that $\imagic(K_n)\le (3+o(1))\, n^2$.  Here we improve on it.

\begin{theorem}\label{th:Wood}
 \begin{equation}\label{eq:imagic}
 \imagic(K_n)\le \left(\frac{288}{121}+o(1)\right)
n^2=(2.380...+o(1))\, n^2.
 \end{equation}
 \end{theorem}
 \smallskip{\it Proof.}  
 Choose $m=\ceil{(\frac{12}{11}+\delta)n}$ for some small constant
$\delta>0$. Take a Sidon set
 \begin{equation}\label{eq:A} 
 \mbox{$A=\{a_1,\dots,a_m\}$ with $1\le a_1<a_2<\cdots<a_m\le
(1+o(1))\, m^2$,}
 \end{equation}
 that is, asymptotically maximum. Explicit such sets were constructed
by Singer~\cite{singer:38} and by Bose and
Chowla~\cite{bose+chowla:62} (Theorems~1 and~3 of
Chapter~II in~\cite{halberstam+roth:s}). 

The case $m=1$ of our Lemma~\ref{lm:EF+L} (or Lemma~1 in Erd\H os and
Freud~\cite{erdos+freud:91}) shows that $A$ is almost uniformly
distributed in $[a_m]$. This implies that if define $T$ to consist of
all triple sums $a_f+a_g+a_h$, $1\le f\le g\le h\le m$, counted with
their multiplicities, then we know the asymptotic distribution of
$T$. We are interested in the interval $[2m^2,3m^2]$, where the
`density' of $T$ at $xm^2$, $2\le x\le 3$, is
 $$
 \int_{x-2}^1 \mathrm{d}y \int_{x-y-1}^1 \mathrm{d}z+o(1)=
\frac{(3-x)^2}2+o(1).
 $$ 
 For example, the number of elements of $T$ lying between $2m^2$ and
$3m^2$ is
 $$
 (1+o(1))\,\binom m3 \int_2^3 \frac{(3-x)^2}2\, \mathrm{d}x =
\left(\frac1{36}+o(1)\right)m^3.
 $$

The interval $I:=[2a_m,(2+\delta)m^2]$ has about $\frac\delta2\,
\binom m3$ elements of $T$, so some $s\in I$ has multiplicity $k\le
(\frac1{12}+o(1))\, m$. For each of the $k$ representations
$s=a_f+a_g+a_h$ remove one of the summands from $A$. Let $B\subset A$
be the remaining set. By removing further elements we can assume that
$|B|=n$.

Label vertices of $K_n$ by the elements of $B$. We want $s$ to be the
magic sum. This determines uniquely the edge labels which are positive
(because $s\ge 2a_m$) and pairwise distinct (because $B\subset A$ is a
Sidon set). Also, as $s\not\in B+B+B$, no edge label equals a vertex
label. As $\delta$ can be chosen arbitrarily small, we obtain
$s=(2+o(1))\,m^2=(\frac{288}{121}+o(1))\, n^2$, proving the
theorem.\qed \medskip

\section{Lower Bound on $\magic(n)$}\label{lower}

For $A\subset \I Z$ let $A\oplus A:=\{a+b\mid a,b\in A,\
a\not=b\}$. We have $A\oplus A\subset A+ A$.

\begin{lemma}\label{lm:sum} Suppose that there is a set $A:=\{a_1=1<a_2<\cdots<a_n\}$ of
integers such that $A\oplus A$ contains an interval of length $m$
(that is, $A\oplus A\supset [k,k+m-1]$ for some $k$). If $a_n\le m$,
then $\magic(n)\ge m-n$.\end{lemma}
 \smallskip{\it Proof.} 
 We will construct an edge-magic graph $G$ on $[n]$ with $m-n$
edges. Label $i\in[n]$ by $l(i):=a_i$. The magic sum will be
$s:=k+m$. For every $a\in A\oplus A$ with $s-a\in [m]\setminus A$
choose a representation $l(i)+l(j)=a$, $1\le i<j\le n$, and add the
pair $\{i,j\}$ (with label $s-a$) to $E(G)$.

Clearly, no two labels are the same. We have 
 $$
 \{s-a\mid a\in A\oplus A\}\supset [m]\supset A
 $$
 So the label set is $[m]$ and we do have an edge-magic graph. The
number of edges is $|[m]\setminus A|=m-n$, as required.\qed \medskip

Mrose~\cite{mrose:79} constructed a set $A\subset[0,10t^2+8t]$ of size
$7t+3$ such that $A+A\supset L:=[0,14t^2+10t-1]$. In fact,
$A=\cup_{i=1}^5 A_i$ is the union of five disjoint arithmetic
progressions. Namely, let
 $$
 [a,(d),b]:=\{a+id\mid i=0,1,\dots,\floor{(b-a)/d}\};
 $$
 then
 \begin{eqnarray*}
 A_1&:=&[0,(1),t],\\
 A_2&:=&[2t,(t),3t^2+t],\\
 A_3&:=&[3t^2+2t,(t+1),4t^2+2t-1],\\
 A_4&:=&[6t^2+4t,(1),6t^2+5t],\\
 A_5&:=&[10t^2+7t,(1),10t^2+8t],
 \end{eqnarray*}
 Fried~\cite{fried:88} independently discovered a similar
construction, giving almost the same bounds.

For any arithmetic progression $B$ we have $|(B+B)\setminus (B\oplus
B)|\le 2$ (because $2b_i=b_{i-1}+b_{i+1}$). Hence, $A\oplus A$
contains all but at most $10$ elements from
$I:=[0,14t^2+10t-1]$. Inspecting each of the ten suspicious elements,
we see that $I\setminus (A\oplus A)=\{0,8t^2+4t-2\}$. Applying
Lemma~\ref{lm:sum} to, for example, the set $\{a+1\mid a\in
A\}\cup\{8t^2+4t-3\}$ with $n=7t+4$, $k=3$, $m=14t^2+10t-1$, we obtain
that $\magic(7t+4)\ge 14t^2+3t-5$ for any $t\ge 1$. Now, the lower
bound in~\eqref{eq:magic} follows from the following lemma.

\begin{lemma}\label{lm:Monotonicity} For any $n$ we have $\magic(n)\le \magic(n+1)$.\end{lemma}
 \smallskip{\it Proof.} 
 Let $G$ be a maximum edge-magic graph of order $n$ with a labelling
$l$.  The graph $G'$ obtained by adding an extra isolated vertex $x$
to $G$ is edge-magic: extend $l$ to $G'$ by defining
$l(x)=v(G)+e(G)+1$.\qed \medskip 

\begin{problem}\label{pr:MLimit} Does the ratio $\magic(n)/n^2$ tend to a limit as
$n\to\infty$?\end{problem}

\section{The Number of Pairwise Sums}\label{PairwiseSums}

The following result is proved via the modification of the argument in
Moser, Pounder and Riddell~\cite[Lemma~1]{moser+pounder+riddell:69}
which in turn is built upon the generating function method of
Moser~\cite{moser:60}. We also refer the reader to a few related
papers: Klotz~\cite{klotz:69}, Green~\cite{green:01}, Cilleruelo,
Ruzsa and Trujillo~\cite{cilleruelo+ruzsa+trjillo:02}, Martin and
O'Bryant~\cite{martin+obryant:02u}.

\begin{theorem}\label{th:1/2} Let $\lambda=\frac14 (2\sqrt2 - 4 + \pi (4-\sqrt2))=
0.323...$\,. Let $n$ be large, $A\subset\I Z$, $m:=|A\setminus[n]|$, and
$k:=|A\cap[n]|$. If $k\ge \lambda m$, then
 \begin{equation}\label{eq:1/2Stronger}
 |(A+A)\cap [2n]|\le n+\frac{|A|^2}4
-\frac{(|A|-\pi m)^2}{(\pi+2)^2}  + o(n),
 \end{equation}
 where the $o(n)$ term depends on $n$ only.
 \end{theorem}
 \smallskip{\it Proof.} 
  Assume that $|A|=O(n^{1/2})$ for otherwise we are done. Let
$A=\{a_1,\cdots,a_{k+m}\}$ with $a_1,\dots,a_k\in[n]$. Correspond to
$A$ its generating function
 $$
 f(x):=\sum_{j=1}^{k+m} x^{a_j}.
 $$
 Let $g(x)=(f^2(x)+f(x^2))/2$. Clearly, the coefficient at $x^j$ in
$g(x)$ is the number of representations of $j$ of the form $a_s+a_t$
with $1\le s\le t\le k+m$.

Let $h(x):=\sum_{j=1}^{2n}
x^j$. Define $\delta_j$, $j\in\I Z$, by the formal identity
 $$
 \sum_{j\in\I Z} \delta_j x^j:=g(x)-h(x).
 $$
 We have $\sum_{j=0}^{2n} \delta_j=\binom {k+m+1}2-2n$. 

Let $t\in[2n-1]$. Then $h(\me^{\pi it/n})=0$, where $i$ is a square
root of $-1$.  Hence,
 $$
 \sum_{j\in\I Z} \delta_j \me^{\pi itj/n}= g(\me^{\pi it/n}).
 $$
 Also observe that each $\delta_j$ is non-negative with the exception
of $j$ lying in $L:=[2n]\setminus(A+A)$ when $\delta_j=-1$. Let
$l:=|L|$.

 Putting all together we obtain, for $t\in[2n-1]$,
 \begin{eqnarray}
 \lefteqn{\frac12\,\left( |f^2(\me^{\pi it/n})|-|f(\me^{2\pi it/n})|\right)
\le |g(\me^{\pi it/n})|\le\Big|\sum_{j\in\I Z\setminus L}
\delta_j\Big| + \Big| \sum_{j\in L} \me^{\pi itj/n}\Big|}
\nonumber\hspace{50pt}\\
 &\le& \sum_{j\in\I Z} \delta_j + 2l + o(n)
 \ =\ \binom{k+m+1}2-2n + 2l + o(n).\label{eq:SumOfDeltas}
 \end{eqnarray}
 Let $z$ denote the right-hand side of~\eqref{eq:SumOfDeltas}, including the
$o(p)$-term.

Let $b_t:=\frac2{t^2-1}$ for even $t>0$ and $b_t:=0$ otherwise.
Clearly, $|f(\me^{2\pi it/n})|\le k+m$ while
 \begin{equation}\label{eq:ModF}
 |f^2(\me^{\pi it/n})|=|f(\me^{\pi it/n})|^2= \Big(\sum_{j\in A} \sin(\pi
ta_j/n)\Big)^2 + \Big(\sum_{j\in A} \cos(\pi ta_j/n)\Big)^2.
 \end{equation}
 Hence, from~\eqref{eq:SumOfDeltas} and~\eqref{eq:ModF} we deduce that
 \begin{eqnarray}
 \frac{\pi}2\,\textstyle (2z)^{1/2}&\ge& \frac\pi2\,\sum_{j\in A} \sin(\pi
a_j/n),\label{eq:111}\\
 b_t\, \textstyle(2z)^{1/2} &\ge&
b_t\,\sum_{j\in A} \cos(\pi ta_j/n),\quad t\in[2,2n-1].\label{eq:211}
 \end{eqnarray}

Note that $\sum_{t=2}^{2n-1} b_t=1-\frac1{2n-1}<1$. By adding~\eqref{eq:111}
and~\eqref{eq:211} we obtain
 \begin{equation}\label{eq:Fourier}
 \left(\frac\pi2+1\right) {\textstyle
(2z)^{1/2}} \ge \sum_{j\in A} \left(\frac\pi2 \sin(\pi
a_j/n) + \sum_{t=2}^{2n-1} b_t \cos(\pi ta_j/n)\right)
 \end{equation}
 It is routine to see that the series $S(x):=\frac\pi2 \sin(x)+
\sum_{t=2}^\infty b_t \cos(tx)$ is the Fourier series of the function
 $$
 r(x)=\left\{\begin{array}{ll}
 1,& 0\le x\le \pi,\\
 1+\pi\sin(x),& \pi\le x\le 2\pi.
 \end{array}\right.
 $$
 (This series appears in~\cite[p.~400]{moser+pounder+riddell:69}.) As
the sum $\sum_{t=2}^\infty |b_t|$ converges and $r(x):\I R/2\pi\I Z\to
\I R$ is a continuous function, it follows from
K\"orner~\cite[Theorem~9.1]{korner:fa} that $S(x)$ converges uniformly
to $r(x)$. Noting that $0\le \pi a_j/n\le \pi$ for any $j\in[k]$, we
conclude that
 \begin{equation}\label{eq:z}
 \left(\frac\pi2+1\right) {\textstyle (2z)^{1/2}}\ge
k+(1-\pi)m+o(m+k).
 \end{equation} 

Assume that $(\pi-1)m> k$ for otherwise we obtain the required by
squaring~\eqref{eq:z}. 

Now,~\eqref{eq:z} is vacuous but we can use the
obvious upper bounds on $|(A+A)\cap [2n]|$ such as $2n$ and
$\binom{k+m+1}2-\frac{m^2}4$. (The latter follows from the fact that
the pairwise sums in $\{x\in A: x>n\}$ and in $\{x\in A\mid x<1\}$ lie
outside $[2n]$.) If neither of these bounds implies~\eqref{eq:1/2Stronger},
then 
 $$
 \frac{(k+m)^2}4- \frac{(k+m(1-\pi))^2}{(\pi+2)^2} < n <
 \frac{(k+m)^2}4 - \frac{m^2}4 + \frac{(k+m(1-\pi))^2}{(\pi+2)^2}.
 $$
 Solving the obtained quadratic inequality in $k$ (and using
$k<m(\pi-1)$), we obtain $k< \lambda m$, as required.\qed \medskip

Note that Theorem~\ref{th:SUpper} easily follows from~\eqref{eq:1/2Stronger}.

\section{Upper Bound on $\magic(n)$}\label{upper}

To prove an upper bound on $\magic(n)$ we study the following function
first. Let $b(k)$ be the largest $n$ such that for some $k$-set
$A\subset \I Z$ we have
 \begin{equation}\label{eq:WeakerAssumption}
 |(A+A)\cap [n]|=(1-o(1))\,n.
 \end{equation}
 It is not hard to see that $b(k)$ has order $\Theta(k^2)$. To state
it formally, we consider the following constant:
 \begin{equation}\label{eq:BSup}
 b_{\sup}:=\limsup_{\e\to0\atop k\to\infty}\, \frac{\max \{n\mid
\exists A\in\binom{\I Z}k,\ |(A+A)\cap [n]|\ge (1-\e)n\}}{k^2}.
 \end{equation}

This definition is related to the question of
Rohrbach~\cite{rohrbach:37} which (when correspondingly reformulated)
asks about $b'(k)$, the largest $n$ such that $[0,n]\subset A+A$ for
some $k$-set $A\subset \I Z_{\ge 0}$. (Note that here $A$ must consist
of \emph{non-negative} integers.) The currently best known upper bound
 $$
 b'(k)\le (0.480...+o(1))\, k^2,
 $$
 is due to Klotz~\cite{klotz:69}. In fact, Klotz's argument gives the
same bound if we weaken the assumption $[0,n]\subset A+A$
to~\eqref{eq:WeakerAssumption}.  The two-side restricted function $b''(k)$
(when we require that $A\subset [0,(\frac12+o(1))n]$) has also been
studied with the present record
 $$
 b''(k)\le (0.424...+o(1))\, k^2,
 $$
 belonging to Moser, Pounder and
Riddell~\cite{moser+pounder+riddell:69} (valid with the weaker
assumption~\eqref{eq:WeakerAssumption} as well).

However, it seems that nobody has considered $b(k)$. Here we fill this
gap as this is the function needed for our application.

\begin{theorem}\label{th:BSup} 
 $$
 b_{\sup}\le \frac12-\frac{2}{(2+(1+2\sqrt2)\, \pi)^2}=\BSupUpper\,.
 $$
 \end{theorem}
 \smallskip{\it Proof.} 
 Let $A\subset \I Z$ have size $k$ and satisfy~\eqref{eq:WeakerAssumption}.
We can assume that $n$ is even. Let $m:=|A\setminus [n/2]|$. As
at least $2\binom{m/2}2=(\frac14+o(1))\, m^2$ sums in $A+A$ fall outside
$[n]$, we have
 \begin{equation}\label{eq:temp11}
 n\le \binom{k}2-\frac{m^2}4+o(k^2)
 \end{equation}
 If $m\ge k/\pi$, then we have
 \begin{equation}\label{eq:n1}
 n\le\left(\frac12-\frac1{4\pi^2}+o(1)\right) k^2 = (0.474...+o(1))\, k^2,
 \end{equation}
 and we are done. Otherwise, by~\eqref{eq:1/2Stronger} we obtain
 \begin{equation}\label{eq:temp21}
 n\le \frac{n}2+\frac{k^2}4-\frac{(k-\pi m)^2}{(\pi+2)^2}+o(k^2).
 \end{equation}
 We conclude that
 $$
 b_{\sup}\le \min_{m\in[0,k]}\left(\frac12-\frac{(m/k)^2}{4},\;
\frac{1}2-\frac{2(1-\pi m/k)^2}{(\pi+2)^2}\right),
 $$
 and the claim routinely follows.\qed \medskip

Let us return to the original problem. Let $l$ be an edge-magic
labelling with the magic sum $s$ of a graph $G$ of order $n$ and size
$m$. Let $A:=l(V(G))$. We have
 \begin{equation}\label{eq:AdditiveBase}
 (A+A)\cap[s-m-n,s-1]\supset \{s-l(\{x,y\})\mid xy\in E(G)\},
 \end{equation}
 that is, $A+A$ contains almost whole interval of length $m+n$
(assuming, obviously, $n=o(m)$). We conclude that $m\le
(b_{\sup}+o(1))\, n^2$, which establishes the upper bound
in~\eqref{eq:magic}.

\section{Asymptotically Maximum Sidon Sequences}

As we have already mentioned the maximum size of a Sidon subset of
$[n]$ is $(1+o(1))\, n^{1/2}$. Erd\H os and
Freud~\cite[Lemma~1]{erdos+freud:91} showed that a set achieving this
bound is almost uniformly distributed among subintervals of $[n]$.
Lindstr\"om~\cite[Theorem~1]{lindstrom:98} proved the analogue of this
result with respect to residue classes.

Here we prove a common generalisation of these results which we will
need in Section~\ref{SLower}. Our proof is based on the method of
Erd\H os and Freud~\cite[Lemma~1]{erdos+freud:91}.

\begin{lemma}\label{lm:EF+L} Let $n$ be large. Let $A$ be an asymptotically maximum
Sidon subset of $[n]$ (that is, having size $(1+o(1))\, n^{1/2}$).
Then for any subinterval $I\subset [n]$ and for any integers $m$ and
$j$, we have
 \begin{equation}\label{eq:EF+L}
 |A\cap I\cap M_j|=\frac{|I|}{m\,n^{1/2}} +o(n^{1/2}).
 \end{equation}
 where $M_j:=\{x\in \I Z\mid x\equiv j\pmod m\}$.
 \end{lemma}
 \smallskip{\it Proof.}  It is enough to prove the lemma for $I=[k]$, an initial
interval, as any other interval is the set-theoretic difference of two
such intervals. Assume that $k=\Omega(n)$ and $m=O(1)$ for
otherwise~\eqref{eq:EF+L} trivially holds.

Choose an integer $t=\Theta(n^{3/4})$. Let $J=\{jm\mid j\in[t]\}$.
For $i\in[-mt+1,n-1]$ let $A_i:=A\cap (I+i)$ and $a_i:=|A_i|$. By the
Sidon property of $A$, the difference set $(A_i-A_i)\cap\I
Z_{>0}\subset J$ has $\binom{a_i}2$ elements; also, a difference
$jm\in J$ is counted $t-j$ times. Hence, we conclude that
 \begin{equation}\label{eq:EF+L1}
 \sum_{j=1}^t (t-j)= \binom t2 \ge \sum_{i=-mt+1}^{n-1}
\binom{a_i}2=\frac12 \sum_{i=-mt+1}^{n-1} a_i^2 - \frac12
\sum_{i=-mt+1}^{n-1} a_i
 \end{equation}

The left-hand size of~\eqref{eq:EF+L1} has magnitude
$t^2=\Theta(n^{3/2})$. All $o(n^{3/2})$-expressions will be dumped
into the error term. In particular, $\sum_i a_i=t|A|$ goes there.

To estimate $\sum_i a_i^2$ we split the summation interval into
smaller parts
 $$
 \mbox{$R_j:=[-mt+1,k]\cap M_j$ and $S_j:=[k+1,n-1]\cap M_j$,\ \
$j\in[m]$.}
 $$
 Now we apply the arithmetic-geometric mean inequality.
 \begin{eqnarray*}
 \sum_{i=-mt+1}^{n-1} a_i^2
 &\ge& \sum_{j\in[m]} \left( \frac{\left(\sum_{i\in R_j}
a_i\right)^2}{|R_j|} + \frac{\left(\sum_{i\in S_j}
a_i\right)^2}{|S_j|}\right)\\
 &=& mt^2 \left(\sum_{j\in[m]}\frac{|A\cap I\cap M_j|^2}k+
\sum_{j\in[m]}\frac{|(A\setminus I)\cap M_j|^2}{n-k}\right)
+o(n^{3/2}).
 \end{eqnarray*}
 (Note that $|R_j|=\frac{k}{m}+O(t)$, $|S_j|=\frac{n-k}m+O(1)$, and
$a_i=O(t^{1/2})$.)

We can estimate the first summand as follows, by using the
arithmetic-geometric mean inequality.
 $$
  \frac{mt^2}{k}\sum_{j\in[m]}|A\cap  I\cap M_j|^2
 \ge \frac{t^2}k \left( \sum_{j\in[m]} |A\cap  I\cap
M_j|\right)^2=\frac{t^2}k |A\cap I|^2.
 $$
 We obtain the analogous bounds for $A\setminus I$. Let
$|A\cap I|=\alpha n^{1/2}$. Then $|A\setminus I|=(1-\alpha+o(1))\,
n^{1/2}$. In summary, starting with~\eqref{eq:EF+L1}, we obtain
 $$
 \binom t2 \ge \frac{t^2}2\left( \frac{|A\cap  I|^2}{k}+
\frac{|A\setminus  I|^2}{n-k}\right) +o(n^{3/2})= t^2\left(\frac12 +
\frac{(\alpha n-k)^2}{k(n-k)}\right) +o(n^{3/2}).
 $$ 

Thus, up to an error term of $o(n^{3/2})$, we must have equality
throughout. We conclude that $\alpha=k/n+o(1)$ and
$a_i=(\alpha/m+o(1))\,n^{1/2}$, which gives the required.\qed \medskip

\section{Lower Bounds on $s(k,n)$}\label{SLower}

We know that the range of interest is $k=\Theta(n^{1/2})$. We will be
proving lower bounds on the following `scaled' one-parameter version
of $s(k,n)$:
 \begin{equation}\label{eq:DefS}
 s(c):=\liminf_{n\to\infty} \frac{s(\floor{cn^{1/2}},n)}n.
 \end{equation}

Note that in~\eqref{eq:DefS} we could have replaced $\floor{cn^{1/2}}$ by
anything of the form $(c+o(1))\, n^{1/2}$ without affecting the value
of $s(c)$. However, we have to write $\liminf$ as the following
question is open.

\begin{problem}\label{pr:limit} Let $c$ be a fixed positive real. Suppose that $n$ tends
to the infinity and $k=(c+o(1))\, n^{1/2}$. Does the ratio $s(k,n)/n$
tend to a limit?\end{problem}

Our lower bound on $s(c)$, provided by the following lemmata, will be
given by different formulae for different ranges of $c$.

The bound~\eqref{eq:QuasiSidonConstruction} of Erd\H os and
Freud~\cite{erdos+freud:91} implies that
 \begin{equation}\label{eq:1}
 s(c)=\frac{c^2}2,\quad c\le 2/\sqrt3.
 \end{equation}
 Their construction can be generalised to give lower bounds on $s(c)$
for larger $c$.

\begin{lemma}\label{lm:SLower1}
 \begin{equation}\label{eq:SLower1}
 s(c)\ge \left\{\begin{array}{ll}
 -\frac{5c^2}8+\frac92 -\frac{6}{c^2}+\frac8{3c^4},& 2/\sqrt3
\le c\le \sqrt2,\\
 \frac{3c^2}8 -\frac32 + \frac{6}{c^2} -\frac{16}{3c^4},& \sqrt2\le
c\le 2.
 \end{array}\right.
 \end{equation}
 \end{lemma}
 \smallskip{\it Proof.}  Let $\alpha=c^2/4$. Choose an integer $m=(\alpha+o(1))\, n$. Let
$A\subset [m]$ be a Sidon set with $(1+o(1))\, m^{1/2}$ elements. The
main idea (which we borrow from Erd\H os and
Freud~\cite{erdos+freud:91}) is to consider the set $X:=A\cup (n-A)$,
where $n-A:=\{n-a\mid a\in A\}$. It is easy to see that, as $A$ is a
Sidon set, all pairwise sums in $A+(n-A)$ are distinct.

However, the set $A+(n-A)$ might intersect  $A+A$. In order to
control the intersection size we introduce some randomness into the
definition of $X$. In what follows, $\e>0$ is a sufficiently small
constant. Let $s,t$ be two integers chosen uniformly and independently
from between $1$ and $\e^2 n$. We define
 $$
 X:=B\cup C,\quad \mbox{where $B:=s+A$ and $C:=n-t-A$.}
 $$ 

Let us compute the densities in $X+X$ which are well defined because
of Lemma~\ref{lm:EF+L}. For example, if we denote
 $$
 \delta_{B+B}(x):=\frac{|(B+B)\cap I|}{|I|},
 $$
 where $I$ is an interval of integers of length $(\e+o(1)) n$ around
$xn$, then
 $$
 \delta_{B+B}(x)=(\mbox{error term})+\left\{\begin{array}{ll}
 \frac{x}{2\alpha},& 0\le x\le \alpha,\\
 -\frac x{2\alpha}+1,& \alpha\le x\le 2\alpha,\\
 0,& \mbox{otherwise,}
 \end{array}\right.
 $$
 where the error term tends to zero if $\e>0$ is sufficiently small
and $n\ge n_0(\e)$. Similarly,
 $$
 \delta_{B+C}(x)=(\mbox{error term})+\left\{\begin{array}{ll}
 0,& 0\le x\le 1-\alpha,\\
 \frac{x}{\alpha}-\frac1{\alpha}+1,& 1-\alpha\le x\le 1.
 \end{array}\right.
 $$
 As the picture is symmetric with respect $x=1$ (given our scaling),
we do not bother about $x\ge 1$ (or about $C+C$).

Thus when one takes some $v\in[n]$ then the probability that $v\in
B+B$ is approximately $\delta_{B+B}(v/n)$. Indeed, this is equivalent
to $v-2s\in A+A$. The case $m=2$ of Lemma~\ref{lm:EF+L} implies that the
number of odd and even elements of $A+A$ in the vicinity of $v$ is
about the same, so their relative density is
$\delta_{A+A}(v)+o(1)$. The analogous claim about the probability of
$v\in B+C$ is also true. Moreover,
 $$
 \Pr\{v\in (B+B)\cap (B+C)\}=\delta_{B+B}(v/n) \times
\delta_{B+C}(v/n) +o(1),
 $$
 because the event is equivalent to $v-2s\in A+A$ and then,
conditioned on this, to $(v-s-n)+t \in A-A$, which has probability
$\delta_{A-A}(\frac{v-s-n}n)+o(1)=\delta_{B+C}(\frac{v}n)+o(1)$.

Hence, by simple inclusion-exclusion, the expectation of $|X+X|$ is at
least
 \begin{equation}\label{eq:SLower2}
 (2+o(1))\, n \int_0^1
(\delta_{B+B}(x)+\delta_{B+C}(x)-\delta_{B+B}(x)\delta_{B+C}(x))\,
\mathrm{d}x.
 \end{equation}
 (Recall that we use the symmetry around $x=1$.) The points $\alpha$,
$1-\alpha$, and $2\alpha$ partition the $x$-range into intervals on
each of which the function in the integral~\eqref{eq:SLower2} is given by
an explicit polynomial in $x$. We have to be careful with the relative
positions of the dividing points: for $\alpha=1/2$ (that is, for
$c=\sqrt2$), the points $\alpha$ and $1-\alpha$ swap places while
$2\alpha$ disappears from the interval. This is why we have two cases
in the bound~\eqref{eq:SLower1} which is obtained by straightforward
although somewhat lengthy calculations (omitted).

Finally observe that there exist $s$ and $t$ such that $|X+X|$ is at
least its expectation, proving the lemma.\qed \medskip

A construction of Rohrbach~\cite[Satz~2]{rohrbach:37} shows that
 \begin{equation}\label{eq:2}
 s(x)=2,\quad\mbox{ if $x\ge 2\sqrt 2$}.
 \end{equation}
 We can extend it for smaller $x$ in the following way.

\begin{lemma}\label{lm:SLower3} Let $c_0:=7/(2\sqrt3)=2.02...$ and
$c_1:=2\sqrt2=2.82...$ Then
 \begin{equation}\label{eq:SLower3}
 s(c)\ge  \left\{\begin{array}{ll}
 \frac{9c^2}{28},& c\le c_0\\
 -c^2 +7\alpha c +\frac{c}{\alpha} -11 \alpha^2- 2-\frac{1}{4\alpha^2},&
c_0\le c\le c_1,
 \end{array}\right.
 \end{equation}
 where $\alpha=\alpha(c)$ is the linear function with
$\alpha(c_0)=\sqrt3/4$ and $\alpha(c_1)=1/\sqrt2$.
  \end{lemma}
 \smallskip{\it Proof.}  Let $k=(c+o(1))\, n^{1/2}$ and let $l:=(3c/14+o(1))\, n^{1/2}$
for $c\le c_0$ and $l:=(\alpha+o(1))\, n^{1/2}$ otherwise.

Let $A:=[l]$, $B:=[n-l+1,n]$. Let $C$ and $D$ be two arithmetic
progressions each of length $\frac{k}2-l$ starting at $(1/2+o(1))\, n$
but with differences $-l$ and $l+1$ respectively. Let $X:=A\cup B\cup
C\cup D$.

All pairwise sums in $A+(C\cup D)$ are distinct, lying within an
interval $[a_0,a_1]$, where $a_0=\frac n2-m+o(n)$ and $a_1=\frac
n2+m+o(n)$, where $m:=(\frac{k}2-l)l$.

Now let us consider $C+D$. Suppose that $c'+d'=c''+d''$ for some
$c'<c''$ in $C$ and $d'>d''$ in $D$. Now, the difference
$c''-c'=d'-d''$ is divisible by both $l$ and $l+1$, hence, it is at
least $l(l+1)$. It is routine to check that $2l^2> m+o(1)\ge l^2$ for
$0< c\le c_1$. This implies that $o(n)$ elements of $C+D$ have
multiplicity at least $3$ and $(\frac k2-2l)^2+o(1)$ elements have
multiplicity $2$ (and all others have multiplicity $1$).

Observe also that $C+D\subset [b_0,b_1]$, where $b_0=n-m+o(n)$ and
$b_1=n+m+o(n)$. 

Let $c\le c_0$. Then $b_0\ge a_1+o(n)$, that is, $A+(C\cup D)$ and
$C+D$ have $o(n)$ elements in common. Therefore, by a sort of symmetry
around $n$, we obtain
 \begin{equation}\label{eq:SLower31}
 |X+X|= 4(k/2-l)l + (k/2-l)^2 - (k/2-2l)^2,
 \end{equation}
 giving the claimed bound.

However, for $c_0\le c\le c_1$, we have $b_0\le a_1+o(n)$. Hence, we
have to subtract from the bound~\eqref{eq:SLower31} twice (by the symmetry)
the number of elements of $C+D$ lying in $[b_0,a_1]$. This correction
term is
 $$
 2\times n\, \int_{b_0/n}^{a_1/n}
\left(\frac{x}{\alpha^2}+\frac{c}{2\alpha} - 1 -\frac{1}{\alpha^2}
\right)\mathrm{d}x +o(n)
 $$
 Computing the value of the integral and plugging it
into~\eqref{eq:SLower31}, the reader should be able to derive the stated
bound.\qed \medskip

\smallskip\noindent{\bf Remark.}  The choice of $l$ for $c_0\le c\le c_1$ in Lemma~\ref{lm:SLower3} is
not best possible. It seems that there is no closed expression for the
optimal choice. So we took a linear interpolation, given the optimal
values for $c=c_0$ and $c=c_1$.\smallskip

Figure~\ref{fg:s} (drawn in \textit{Mathematica}) contains the
graphical summary of our findings.

\section{Differences}\label{differences}

Similar questions can be asked about differences. For example, let us
define
 $$
 d(k,n):=\max\left\{ |A-A|\mid A\in {[n]\choose k}\right\}.
 $$

The obvious upper bounds are $2n-1$ and $k(k-1)+1$ (where the last
summand $1$ counts $0\in A+A$). These bounds can be improved when
$\sqrt n\le (1+o(1))k\le \frac32\, \sqrt n$ as the following theorem
demonstrates. 

\begin{theorem}\label{th:DifferenceUpper} Let $n$ be large and $k\ge \sqrt n$. Then
 \begin{equation}\label{eq:DifferenceUpper}
 d(k,n)\le 2k\sqrt n - n +o(n).
 \end{equation}
 \end{theorem}
 \smallskip{\it Proof.} 
 Let $c:=k/\sqrt n>1$. Assume that $c-1=\Theta(1)$ for otherwise we are
trivially done. Define $t:=\floor{(c-1)n}$,
 $$
 A_i:=A\cap [i,i+t-1],\mbox{ and } a_i:=|A_i|,\quad i\in[2-t,n].
 $$ 
 
Let $\C X$ consist of all quadruples $(a,b,i,x)$ such that $x=a-b>0$
and $a,b\in A_i$. Using the identity $\sum_{i=2-t}^{n} a_i=kt$ and the
quadratic-arithmetic mean inequality, we obtain
 \begin{equation}\label{eq:BigA}
 |\C X|=\sum_{i=2-t}^{n} \binom{a_i}2= \frac12\, \sum_{i=2-t}^{n}
a_i^2\ - \frac{kt}2 \ge
 (1+o(1)) \frac{(kt)^2}{2(n+t)}.
 \end{equation}

For $x\in \I N$, let $g_x$ be the number of representations $x=a-b$
with $a,b\in A$. Then, each $x\in [t-1]$ is included in $g_x(t-x)$
quadruples. Hence,
 \begin{equation}\label{eq:BigA1}
 |\C X|\le \sum_{x=0}^{t-1} (t-x)g_x.
 \end{equation}

The above sum can be bounded by
$\sum_{i=0}^{t-i}(t-i)=(\frac12+o(1))\, t^2$
plus $\frac12\, t(k^2-|A-A|)$. Putting all together we obtain:
 $$
 \frac{(kt)^2}{2(n+t)}\le \frac{t^2}2 + \frac{t(k^2-|A-A|)}2+o(n^2).
 $$

Routine simplifications yield the claim.\qed \medskip

Let us briefly discuss the lower bounds on
 $$
 d(c):=\liminf_{n\to\infty} \frac{d(\floor{cn^{1/2}},n)}n.
 $$

Sidon sets show that $d(c)=c^2$ for $0\le c\le1$. 

\begin{lemma}\label{lm:d1} For $1\le c\le \sqrt2$,
 $$
 d(c)\ge  - \frac{c^4}3 + 2c^2 - 2+ \frac{4}{3c^2}.
 $$
 \end{lemma}
 \smallskip{\it Proof.} 
  Let $\beta=1/c^2$ and $b=\floor{\beta n}$. Let $B\subset [b]$ be a
maximal Sidon set. Let $C=[n]\cap (B+b)$ and $A=B\cup (C+t)$,
where $t$ is a small random integer. As $B$ is uniformly distributed
in $[b]$, it is easy to see that $|A|=(c+o(1))\, \sqrt n$ is as
required.

All differences in $C-B$ are pairwise distinct. So, the densities of
$B-B$ and $C-B$ at $xn$, $0\le x\le 1$, are respectively
$f(x)=1-x/\beta$ if $0\le x\le \beta$ (while $f(x)=0$ for $x\ge
\beta$) and
 $$
 g(x)=\left\{
 \begin{array}{ll}
 x/\beta,& 0\le x\le 1-\beta,\\
 (1-\beta)/\beta,& 1-\beta\le x\le \beta,\\
 (1-x)/\beta, & \beta\le x\le 1,\\
 0,& \mbox{otherwise.}
 \end{array}\right.
 $$
 (Note that $C-C\subset B-B$, so there is no point to consider $C-C$.)

Now, similarly to our analysis in Lemma~\ref{lm:SLower1}, the expected
size of $A-A$ is
 $$
 (2+o(1)) n \int_0^1 (f(x)+g(x)-f(x)g(x))\md x = n\,
\Big(\frac{4\beta}3-2+\frac 2\beta- \frac1{3\beta^2}+o(1)\Big).
 $$
 By taking $t$ so that $|A-A|$ is at least its expectation, we complete the
 proof.\qed \medskip

The following construction provides best known lower bounds for the
remaining values of $c$.

Choose some $\beta\le c$ (to be specified later). Let $b=\floor{\beta
\sqrt n}$. Define $B=[b]$. Let $C$ and $D$ be arithmetic progressions
of length $\frac{c-\beta}2 \sqrt n$ starting at
$(1-\frac{\beta(c-\beta)}2+o(1)) n$ but the differences $-b$ and
$b-1$ respectively. (Thus, for example, $D$ ends around $n$.) Let
$A=B\cup C\cup D$. Clearly, $(C\cup D)-B$ covers an interval
$[(1-\beta(c-\beta)+o(1))n,n-1]$. Also, the distribution of $D-C$ can be
explicitely written, which allows us to compute $|A-A|$
asymptotically.

For $c\ge 2$, we can ensure that $A-A = [-n+1,n-1]$; thus $d(c)=2$
then. For $\sqrt2\le c\le 3/2$, the optimal choice is $\beta=c/3$,
giving
 $$
 d(c)\ge  2c^2/3,\quad \sqrt2\le c\le 3/2.
 $$

\begin{figure}[h]
\begin{minipage}[t]{6.0cm}
\begin{center}
\scalebox{0.59}[0.8]{\includegraphics{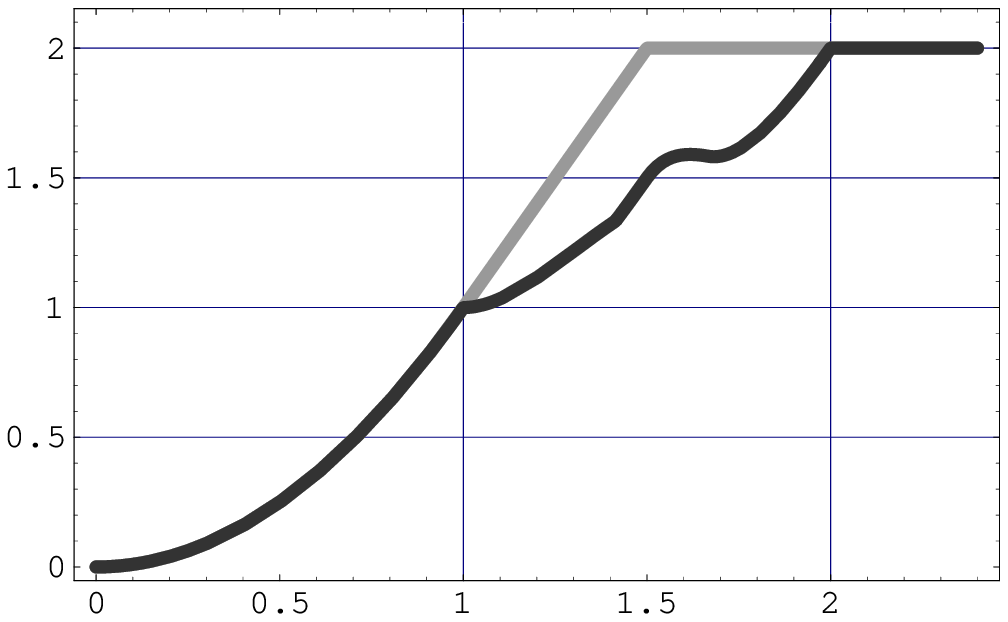}}
\end{center}
\end{minipage}
\hfill
\begin{minipage}[t]{6.0cm}
\begin{center}
\scalebox{0.59}[0.8]{\includegraphics{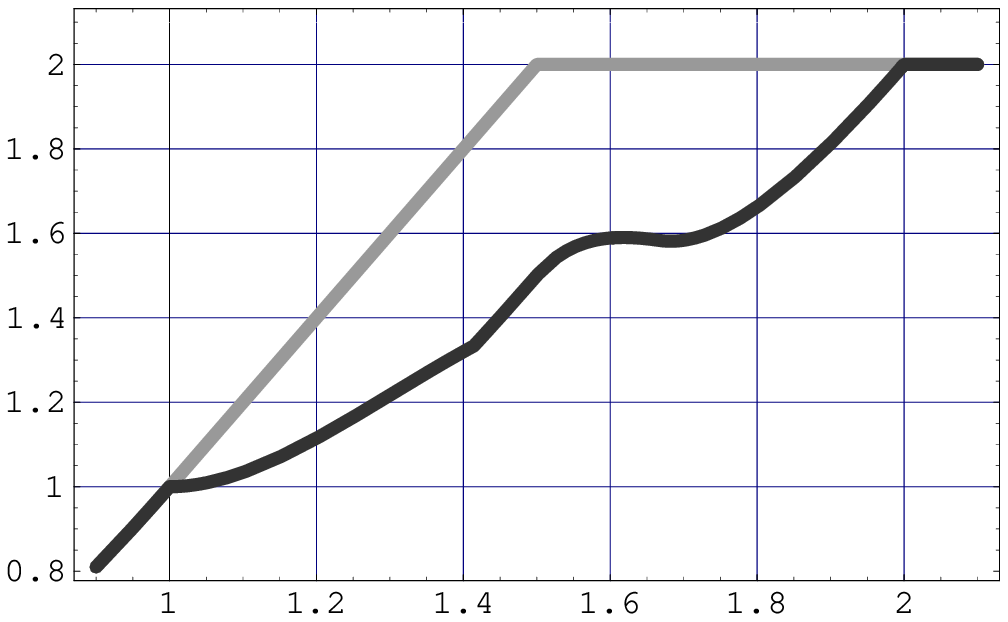}}
\end{center}
\end{minipage}
\caption{\small Our bounds on $d(k,n)$: $x=kn^{-1/2}$, 
 $y=d(k,n)/n$.}\label{fg:d}
\end{figure}

Unfortunately, it seems that there is no closed formula the optimal
$\beta=\beta(c)$ for other values of $c$. (And, in fact, $\beta(c)$ is
not a continuous function.) But, as an illustration, we choose
$\beta(c)=c-1$, the linear interpolation given the optimal choices
$\beta(3/2)=1/2$ and $\beta(2)=1$. Routine calculations give us the
following lower bounds.
 $$
 d(c)\ge \left\{\begin{array}{ll}
 \frac{4c^3-19c^2+34c-21}{2(c-1)^2},& \frac32\le c\le \frac53,\\
 \frac{2c^3-5c^2+2c+2}{(c-1)^2},& \frac53\le c\le 2.
 \end{array}\right.
 $$

Figure~\ref{fg:d} contains the graphs of our bounds.

\begin{problem}\label{pr:DiffBasis} Compute $d(n)$, the smallest size of $A\subset [n]$
such that $A-A\supset [-n+1,n-1]$. The same question about $d'(n)$
when we require that $|A-A|=(2+o(1))n$ only. Is $d(n)=(1+o(1))\,
d'(n)$?\end{problem}

At the moment we know only that $d'(n)\le d(n)$ lie between $\frac
32\, n$ and $2n$.

\begin{problem}\label{pr:DiffLimit} Does the ratio $d(k,n)/n$ tend to a limit as
$n\to\infty$ and $k=(c+o(1))\, n^{1/2}$ where $c$ is fixed?\end{problem}

\end{document}